\documentclass[12pt]{article}
\usepackage{amsmath,amsfonts}
\begin{document}

\baselineskip 20pt
\title{Fatou's interpolation theorem implies the Rudin-Carleson theorem
}

\author{Arthur A.~Danielyan}

\maketitle

\begin{abstract}

\noindent The purpose of this paper is to show that the Rudin-Carleson interpolation theorem is a direct corollary of Fatou's much older interpolation theorem (of 1906).

\end{abstract}

\begin{section}{Introduction.}

Denote by $\Delta$ and $T$ the open unit disk and the unit circle in the complex plane,
respectively. Recall that the disk algebra $A$ is the algebra of all
functions on the closed unit disk $\overline \Delta$ that are
analytic on $\Delta$.
 
The following theorem is fundamental; in particular it implies the F. and M. Riesz theorem on analytic measures (cf. \cite{Koos}, pp. 28-31). 

\vspace{0.25 cm}

{\bf Theorem A (P. Fatou, 1906).}\ {\it Let $E$ be a closed subset of $T$ such that $m(E)=0$ ($m$ is the Lebesgue measure on $T$). 
Then there exists a function $\lambda_E(z)$ in the disk algebra $A$ such that $\lambda_E(z)=1$ on $F$ and $|\lambda_E(z)|<1$ on $T \setminus E$.}

\vspace{0.25 cm}

In its original form Fatou's theorem states the existence of an element of $A$ which 
vanishes precisely on $F$, but it is equivalent to the above version (cf. \cite{Koos}, p. 30, or \cite{Hoff}, pp. 80-81).

The following famous
 theorem, due to W. Rudin \cite{Rud} and L. Carleson \cite{Carl},
 has been the starting point of many investigations in complex and functional analysis (including several complex variables).

\vspace{0.25 cm}

{\bf Theorem B (Rudin - Carleson).}\ {\it Let $E$ be a closed set of measure zero on $T$ and let $f$ be a continuous (complex valued) function on $E$. 
Then there exists a function $g$ in the disk algebra $A$ agreeing with $f$ on $E.$}

\vspace{0.25 cm}


It is obvious from Theorem A and Theorem B that for any $\epsilon>0$ one can choose the extension function $g$ in Theorem B such that it is bounded by   $||f||_E+ \epsilon$,
where  $||f||_E$ is the sup norm of $f$ on $E.$ Rudin has shown that one can even choose the function $g$ such that it is bounded by $||f||_E$.

Quite naturally, as mentioned already by Rudin,
Theorem B may be regarded as a strengthened  form of Theorem A (cf. \cite{Rud}, p. 808). 

The present paper
shows that Theorem B also is an elementary corollary of Theorem A. 
To be more specific, we present a brief proof of Theorem B merely using Theorem A and the Heine - Cantor theorem
(from Calculus 1 course); this approach may find further applications. 
We use a simple argument based on uniform continuity, which has been known (at least since 1930s) in particular to
M.A. Lavrentiev \cite{Lav}, M.V.  Keldysh, and S.N. Mergelyan,
but 
has not been used for the proof of Theorem B before.



\end{section}
\begin{section}{Proof of Theorem B}


Let $\epsilon>0$ be given.  
By uniform continuity we cover  $E$ by disjoint open intervals $I_k \subset T$ of a finite number $n$ such that  $|f(z_1)-f(z_2)| < \epsilon$
 for any $z_1, z_2 \in E \cap I_k$ ($k=1,2,..., n$). Denote $ E_k= E \cap I_k$ and let
  $\lambda_{E_k}(z) $ be the function provided by Theorem A.
 Fix a natural number $N$ so large that $|\lambda_{E_k}(z)|^N<\frac{\epsilon}{n}$ on $T\setminus I_k$ for all $k.$ Fix a point $t_k \in E_k$ for each $k$ and denote
 $h(z)=\sum_{k=1}^nf(t_k) [\lambda_{E_k}(z)]^{N}$. Obviously the function $h \in A$ is bounded on $T$ by the number $ (1+\epsilon) ||f||_E$ and 
$|f(z)-h(z)| < \epsilon(1+  ||f||_E) $ if $z\in E$. Replacing $h$ by $\frac{1}{1+\epsilon}h$ allows to assume that $h$ is bounded on $T$ simply by $||f||_E$ and $|f(z)-h(z)| < \epsilon(1+ 2 ||f||_E) $ if $z\in E$. Letting $\epsilon = \frac{1}{m}$ provides a sequence $\{h_m\}$, $h_m \in A$, which is uniformly bounded on $T$ by $  ||f||_E$ and uniformly
 converges to $f$ on $E$. 
 
 To complete the proof, we use the following known steps (cf. e.g. \cite{Dan}).
Let $\eta>0$ be given and let $\eta_n>0$ be such $\sum \eta_n < \eta$. We can find $H_1=h_{m_1} \in A$ such that  $|H_1(z)| \leq ||f||_E$ on $T$ and $|f(z)-H_1(z)| < \eta_1$ on $E$. 
Letting $f_1=f - H_1$ on $E$, the same reasoning yields $H_2 \in A$ with  $|H_2(z)| \leq ||f_1||_E < \eta_1$ on $T$ and $|f_1(z)-H_2(z)| < \eta_2$ on $E$. 
Similarly we find $H_n \in A$ for $n=3, 4, ...,$ with appropriate properties. The convergence of the series $||f||_E +  \eta_1 + \eta_2 + ... $ implies that the series $\sum H_n(z) $  converges uniformly on $\overline \Delta$ to a function $g \in A$, which is bounded by $||f||_E+\eta$. On $E$ holds $|f-g|=|(f-H_1)-H_2-...-H_n-...|=|(f_1-H_2)-H_3-...-H_n-...|=...=|(f_{n-1} -H_n)-...| \leq \eta_n + \sum_{k=n}^{\infty} \eta_k$. 
Since  $\lim_{n \rightarrow \infty } (\eta_n + \sum_{k=n}^{\infty} \eta_k)=0$, it follows that $g=f$ on $E$, which completes the proof.

\vspace{0.25 cm}

{\bf Remark.} The known proofs of Theorem B 
use Theorem A and a polynomial approximation theorem  (
cf. \cite{Gar}, p. 125; or  
\cite{Hoff}, pp. 81-82).
The latter is needed to approximate  $f$ on $E$ by elements of the disc algebra $A$.
 The above proof uses just Theorem A to provide such approximation of $f$ on $E$ by elements of $A$, which in addition are bounded by $||f||_E$  on $T$.



\end{section}

\begin{minipage}[t]{6.5cm}
Arthur A. Danielyan\\
Department of Mathematics and Statistics\\
University of South Florida\\
Tampa, Florida 33620\\
USA\\
{\small e-mail: adaniely@usf.edu}
\end{minipage}

\end{document}